\documentclass[11pt, a4paper, abstract=yes]{scrartcl}

\usepackage{hyperref}
\usepackage{amsthm,amssymb,amsmath}
\usepackage{hyperref}
\usepackage{graphicx}
\usepackage{iftex}
\ifPDFTeX
\usepackage[T1]{fontenc}\usepackage{lmodern}
\else
\usepackage{fontspec}
\IfFontExistsTF{Latin Modern Roman}{\setmainfont{Latin Modern Roman}}{}
\IfFontExistsTF{Latin Modern Sans}{\setsansfont{Latin Modern Sans}}{}
\IfFontExistsTF{Latin Modern Mono}{\setmonofont{Latin Modern Mono}}{}
\fi
\theoremstyle{definition}
\newtheorem{definition}{Definition}

\theoremstyle{remark}
\newtheorem{remark}[definition]{Remark}
\theoremstyle{plain}
\newtheorem{lemma}[definition]{Lemma}

\newtheorem{proposition}[definition]{Proposition}
\newtheorem{theorem}[definition]{Theorem}
\newtheorem{corollary}[definition]{Corollary}

\author{Jan Snellman\thanks{Department of Mathematics, Linköping University, 581 83 Linköping, Sweden. \texttt{jan.snellman@liu.se}. ORCID: \href{https://orcid.org/0009-0002-6676-5068}{0009-0002-6676-5068}}}
\usepackage{booktabs}
\KOMAoptions{captions=nooneline}
\usepackage{tikz}
\usepackage{marginnote}
\definecolor{leangreen}{RGB}{0,140,60}
\definecolor{leanorange}{HTML}{E65100}
\newcommand{\leanfull}{\tikz[baseline=-0.6ex]\fill[leangreen] (0,0) circle (0.55ex);}
\newcommand{\leannative}{\tikz[baseline=-0.6ex]\fill[leanorange] (0,0) circle (0.55ex);}
\newcommand{\leanfullnative}{\leanfull\,\leannative}
\usetikzlibrary{calc}
\newcommand{\leaneq}[1]{\tikz[remember picture,overlay]{\coordinate (leanhere) at (0,0.6ex);\node[anchor=east,inner sep=0pt] at ($(current page.west |- leanhere)+(1in+\hoffset+\oddsidemargin-\marginparsep,0)$) {#1};}}
\DeclareRobustCommand{\leannone}{\tikz[baseline=-0.6ex]\draw[black!55,line width=0.6pt] (0,0) circle (0.5ex);}
\newcommand{\leanmark}[1]{\leavevmode{\reversemarginpar\marginnote{\raggedleft#1}}}
\makeatletter\@ifpackageloaded{showkeys}{\renewcommand*\showkeyslabelformat[1]{\fbox{\normalfont\tiny\ttfamily#1}\kern1.5em}}{}\makeatother
\newcommand{\Iver}[1]{[\![\,#1\,]\!]}
\newcommand{\SYT}{\operatorname{SYT}}
\newcommand{\Am}{A^{(m)}}
\newcommand{\Bm}{B^{(m)}}
\newcommand{\Cm}{C^{(m)}}
\newcommand{\Gm}{G^{(m)}}
\newcommand{\tCm}{\widetilde{C}^{(m)}}
\date{\today}
\title{The Growth Constant for Multichains in Young's Lattice}
\hypersetup{
 pdfauthor={},
 pdftitle={The Growth Constant for Multichains in Young's Lattice},
 pdfkeywords={},
 pdfsubject={},
 pdfcreator={},
 pdflang={English}}
\usepackage[doi=false,url=false,isbn=false,style=numeric]{biblatex}
\begin{document}

\maketitle
\begin{abstract}
Azam and Richmond showed that, over the partitions \(\lambda\) of \(n\) with exactly \(k\) parts,
the average size of the interval \([\emptyset,\lambda]\) in Young's lattice grows like \(G_k n^k\),
with \(G_k\) a weighted sum over an interval of partitions. We give a short proof that does not use
their recursion, with an error term, and for multichains
\(\lambda \supseteq \mu^1 \supseteq \dots \supseteq \mu^{m-1}\) of any length: the average number
of them below \(\lambda\) is \(\Gm_k n^{(m-1)k}(1+O(1/n))\), where \(\Gm_k\) is a weighted sum over
the standard Young tableaux of the \(m \times k\) rectangle. For \(k=2\) the constant has a closed
form.
\end{abstract}


\medskip
\section{Introduction}
\label{sec:orgee95175}

We write \(\lambda \vdash n\) when \(\lambda\) is a partition of \(n\), \(\ell(\lambda)\) for its
number of parts, \(\mu \subseteq \lambda\) for containment of diagrams, and \((k^m)\) for the
partition with \(m\) parts equal to \(k\), whose diagram is the rectangle with \(m\) rows of
length \(k\). Further \([x^n]f\) is the coefficient of \(x^n\) in \(f\), and \(\Iver{S}\) is the
Iverson bracket, \(1\) if the statement \(S\) holds and \(0\) otherwise. Let \(c_{k,n}\) be the number of partitions of \(n\) with exactly
\(k\) parts. For a partition \(\lambda\) let \([\emptyset,\lambda]\) be the set of partitions \(\mu \subseteq \lambda\), an
interval in Young's lattice. Azam and Richmond \autocite{azam_generating_2023} study
\[
C_{k,n} = \sum_{\lambda \vdash n,\ \ell(\lambda)=k} \#[\emptyset,\lambda],
\qquad A_{k,n} = \frac{C_{k,n}}{c_{k,n}},
\]
the average size of such an interval. Their Theorem 3 states that \(A_{k,n} \sim G_k n^k\) as
\(n \to \infty\), where
\begin{equation}
\label{eq-ar}
G_k = \frac{(k-1)!}{(2k-1)!} \sum_{\lambda \in I_k} \prod_{i=1}^{k} \frac{1}{\lambda_i},
\qquad I_k = [(k,k-1,\dots,1),\,(k^k)].
\end{equation}
They obtain this from a rational recursion for a multivariate generating function, of which
that of the \(C_{k,n}\) is a specialisation; the recursion was the subject of the earlier comments \autocite{snellman_multi_2025}. The purpose of this
note is to point out that the asymptotics follow directly from Stanley's formula for
\(P\)-partitions, and that the argument does not see the number of rows.

Fix \(m \ge 2\). A \emph{multichain} below \(\lambda\) is a sequence
\(\lambda \supseteq \mu^1 \supseteq \dots \supseteq \mu^{m-1}\) of partitions. Put
\[
\Cm_{k,n} = \sum_{\lambda \vdash n,\ \ell(\lambda)=k}
\#\{\lambda \supseteq \mu^1 \supseteq \dots \supseteq \mu^{m-1}\},
\qquad \Am_{k,n} = \frac{\Cm_{k,n}}{c_{k,n}},
\]
so that \(C^{(2)}_{k,n} = C_{k,n}\) and \(A^{(2)}_{k,n} = A_{k,n}\).

Let \(\SYT(k^m)\) be the set of standard Young tableaux of shape \((k^m)\). For \(T \in \SYT(k^m)\) and \(1 \le t \le mk\) let
\[
a_t(T) = \text{the number of entries of the first row of } T \text{ that are} \le t .
\]
Thus \(1 = a_1(T) \le a_2(T) \le \dots \le a_{mk}(T) = k\).

\begin{theorem}
\leanmark{\leanfull}Let \(k \ge 1\) and \(m \ge 2\) be fixed.\footnote{The discs in the margin record which statements have been checked in Lean; see Section \ref{sec-marks}.} As \(n \to \infty\),
\[
\Am_{k,n} = \Gm_k\, n^{(m-1)k}\,\bigl(1 + O(1/n)\bigr),
\qquad \Gm_k = \frac{k!\,(k-1)!}{(mk-1)!}\, \Bm_k,
\]
where
\begin{equation}
\label{eq-Ba}
\Bm_k = \sum_{T \in \SYT(k^m)} \prod_{t=1}^{mk} \frac{1}{a_t(T)} .
\end{equation}
Moreover
\begin{equation}
\label{eq-Bb}
\Bm_k = \frac{1}{k!} \sum_{T' \in \SYT(k^{m-1})} \sum_{s \in S(T')} \prod_{i=1}^{(m-1)k} \frac{1}{s_i},
\end{equation}
where \(S(T')\) is the set of integer sequences \(1 \le s_1 \le \dots \le s_{(m-1)k} \le k\) with
\(s_i \ge a_i(T')\) for all \(i\), where \(a_i(T')\) is, as above, the number of entries \(\le i\) in
the first row of \(T'\).
\label{thm-main}
\end{theorem}

For \(m=2\) there is one tableau \(T'\), a single row, and \(a_i(T') = i\). Then \(s \in S(T')\) exactly
when \((s_k,\dots,s_1) \in I_k\), so \eqref{eq-Bb} gives \(G^{(2)}_k = G_k\) and the theorem contains
that of Azam and Richmond, with an error term (Corollary \ref{cor-ar}). The elements of \(I_k\) are counted by the Catalan numbers, which is the ``weighted sum over
Catalan objects'' of \autocite[Remark 5]{azam_generating_2023}. For \(m=3\) a second Catalan family
appears: the outer sum in \eqref{eq-Bb} is over the tableaux \(T'\) of shape \((k,k)\), and for
each of them the inner sum is over the sequences above a staircase determined by \(T'\).

Table \ref{tab-G} lists some values. In two cases there is a closed form: \(\Gm_1 = 1/(m-1)!\),
which is elementary, and

\begin{proposition}
\leanmark{\leanfull}\(\displaystyle \Bm_2 = \frac{1}{4^m}\binom{2m}{m}\), and hence
\(\displaystyle \Gm_2 = \frac{2}{4^m\,(2m-1)!}\binom{2m}{m}\).
\label{prop-two}
\end{proposition}

\begin{table}[htbp]
\caption{\label{tab-G}\leannone{} The growth constants \(\Gm_k\). The column \(m=2\) agrees with Table 1 of Azam and Richmond, which goes to \(k=7\).}
\centering
\begin{tabular}{cccc}
\toprule
\(k\) & \(G^{(2)}_k\) & \(G^{(3)}_k\) & \(G^{(4)}_k\)\\
\midrule
1 & \(1\) & \(1/2\) & \(1/6\)\\
2 & \(1/8\) & \(1/192\) & \(1/9216\)\\
3 & \(49/6480\) & \(1579/117573120\) & \(77423/8380611993600\)\\
4 & \(1597/5806080\) & \(2072569/158921234841600\) & \(65198099/342720635470243430400\)\\
\bottomrule
\end{tabular}
\end{table}
\subsection{Formal verification}
\label{sec-marks}
Every numbered statement and table of this note carries one mark in the margin, which says
whether it has been formalised in Lean 4 with Mathlib. The marks refer to statements; a proof is marked
only by its heading, which says so when it is the route that the formalisation took.
\begin{itemize}
\item \leanfull{} Formalised and checked, using no axioms beyond the three standard ones. A mark
level with a displayed formula, as for \eqref{eq-diff} and \eqref{eq-gf}, refers to that
formula alone.
\item \leannative{} A finite computation that Lean ran as compiled code, which trusts the Lean
compiler as well. This occurs only in the example of Section \ref{sec-example}: the counts for
\(n=4\), every entry of Table \ref{tab-ex}, the sets \(S(T')\), and the two generating
functions, which follow from \eqref{eq-gf} and \eqref{eq-diff} once the tableaux are listed.
\item \leannone{} Not formalised: Table \ref{tab-G}, apart from its row \(k = 2\), which is Proposition
\ref{prop-two}, and Remark \ref{rem-method}. Also not
formalised are the first proofs of Proposition \ref{prop-two} and of Lemma \ref{lem-coeff}, and the
numerical ratios at the end of the example.
\end{itemize}
In the formalisation a multichain is, by definition, a plane partition as in Section
\ref{sec-proof}, and \(c_{k,n}\) is the case \(m=1\) of the same count, which is proved equal to
the number of partitions of \(n\) with \(k\) parts in the sense of Mathlib.
\section{Proof of the theorem}
\label{sec-proof}
Let \(P = [m] \times [k]\) be the set of cells \((i,j)\), row \(i\) and column \(j\). A \emph{plane
partition} on \(P\) is an array \(\pi\) of nonnegative integers on \(P\), weakly decreasing along
rows and down columns. Its rows \(\lambda = \pi(1,\cdot) \supseteq \mu^1 = \pi(2,\cdot) \supseteq
\dots \supseteq \mu^{m-1} = \pi(m,\cdot)\) form a multichain with \(\ell(\lambda) \le k\), and every
such multichain arises once. Let \(\tCm_{k,n}\) be the number of plane partitions on \(P\) whose
first row has sum \(n\). Since \(\ell(\lambda) \le k-1\) means that the last column of \(\pi\) is
zero,
\begin{equation}
\label{eq-diff}
\leaneq{\leanfull}
\Cm_{k,n} = \tCm_{k,n} - \tCm_{k-1,n}, \qquad \tCm_{0,n} = \Iver{n=0}.
\end{equation}

A tableau \(T \in \SYT(k^m)\) is a bijection \(P \to \{1,\dots,mk\}\) increasing along rows and
down columns. Let \(w_t\) be the cell with entry \(t\), and let \(D(T)\) be the set of \(t < mk\)
such that \(w_{t+1}\) lies in a row strictly above the row of \(w_t\).

\begin{lemma}
\leanmark{\leanfull}With a variable \(x_c\) for each cell \(c\), and \(u_t(T) = x_{w_1} x_{w_2} \cdots x_{w_t}\),
\[
\sum_{\pi} \prod_{c \in P} x_c^{\pi(c)}
= \sum_{T \in \SYT(k^m)} \frac{\prod_{t \in D(T)} u_t(T)}{\prod_{t=1}^{mk} \bigl(1-u_t(T)\bigr)},
\]
the sum on the left being over all plane partitions on \(P\).
\label{lem-stanley}
\end{lemma}

This is Stanley's formula for \(P\)-partitions \autocite[Theorem 3.15.5]{stanley_ec1}, for the poset \(P\) ordered from north-west to south-east
with the natural labelling \(\omega\) below. We include the short proof, since the set \(D(T)\) depends on
conventions.

\begin{proof}
Put \(\omega(i,j) = (i-1)k+j\). If \(c \ne c'\) and \(c\) is weakly north-west of \(c'\), then
\(\omega(c) < \omega(c')\). Given \(\pi\), list the cells by decreasing value of \(\pi\), and cells
of equal value by increasing \(\omega\). A cell weakly north-west of another is listed before it,
so the list \(w_1,\dots,w_{mk}\) is a tableau \(T(\pi)\). Put \(d_t = \pi(w_t)-\pi(w_{t+1})\) for
\(t<mk\) and \(d_{mk} = \pi(w_{mk})\). Then \(d_t \ge 0\), and \(d_t \ge 1\) whenever
\(\omega(w_t) > \omega(w_{t+1})\). Since \(\pi(w_t) = d_t + \dots + d_{mk}\),
\[
\prod_{c} x_c^{\pi(c)} = \prod_{t=1}^{mk} u_t(T)^{d_t}.
\]
Conversely, a tableau \(T\) and integers \(d_t \ge 0\), with \(d_t \ge 1\) whenever
\(\omega(w_t) > \omega(w_{t+1})\), define a plane partition \(\pi(w_t) = d_t + \dots + d_{mk}\) with
\(T(\pi) = T\): along a run of equal values \(d_t = 0\), so \(\omega\) increases, which is how
\(T(\pi)\) lists the run. So \(\pi \mapsto (T(\pi),d)\) is a bijection. Finally
\(\omega(w_t) > \omega(w_{t+1})\) if and only if \(t \in D(T)\): if \(w_{t+1}\) is in the row of
\(w_t\) it is to the right of it, and \(\omega\) increases; if it is in a lower row \(\omega\)
increases; if it is in a higher row \(\omega\) decreases. Summing the geometric series over each \(d_t\)
gives the formula.
\end{proof}

\begin{lemma}
\leanmark{\leanfull}Let \(a_1,\dots,a_N\) be positive integers, at least one of them equal to \(1\), and let
\(e \ge 0\) be an integer. Then, as \(n \to \infty\),
\[
[x^n]\, \frac{x^e}{\prod_{t=1}^{N}(1-x^{a_t})}
= \frac{n^{N-1}}{(N-1)!\; a_1 \cdots a_N} + O(n^{N-2}).
\]
\label{lem-coeff}
\end{lemma}

\begin{proof}
Call the function \(f\). Its poles are roots of unity \(\zeta\), of order
\(N_\zeta = \#\{t : \zeta^{a_t}=1\}\). So \(N_1 = N\), while \(N_\zeta \le N-1\) for
\(\zeta \ne 1\), because some \(a_t\) is \(1\). By partial fractions,
\(f(x) = p(x) + \sum_\zeta \sum_{j=1}^{N_\zeta} c_{\zeta,j}(1-x/\zeta)^{-j}\) with a polynomial
\(p\), and the coefficient of \(x^n\) in \((1-x/\zeta)^{-j}\) is \(\binom{n+j-1}{j-1}\zeta^{-n}\),
which is \(O(n^{j-1})\). Hence \([x^n]f = c_{1,N}\binom{n+N-1}{N-1} + O(n^{N-2})\) with
\[
c_{1,N} = \lim_{x \to 1} (1-x)^N f(x) = \prod_{t=1}^{N} \lim_{x\to 1}\frac{1-x}{1-x^{a_t}}
= \prod_{t=1}^{N}\frac{1}{a_t}.
\]
\end{proof}

The formal proof of Lemma \ref{lem-coeff} is elementary and gives explicit bounds.

\begin{proof}[Second proof of Lemma \ref{lem-coeff} (the formalised route)]
We may assume \(a_N = 1\). For \(0 \le s \le N-1\) let \(c_s(n)\) be the coefficient of \(x^n\) in
\((1-x)^{-1}\prod_{t \le s}(1-x^{a_t})^{-1}\), and put \(P_s = a_1 \cdots a_s\) and
\(S_s = a_1 + \dots + a_s\). We claim that for all \(n \ge 0\)
\[
\frac{n^s}{s!\,P_s} \;\le\; c_s(n) \;\le\; \frac{(n+S_s)^s}{s!\,P_s}.
\]
For \(s=0\) this says \(c_0(n) = 1\). With \(a = a_{s+1}\), multiplying by \(1-x^a\) gives
\(c_{s+1}(n) = c_s(n) + c_{s+1}(n-a)\), the last term being absent when \(n < a\), and the claim
for \(s+1\) follows by induction on \(n\) from the inequalities
\[
y^{s+1} + (s+1)\,a\,y^s \;\le\; (y+a)^{s+1} \;\le\; y^{s+1} + (s+1)\,a\,(y+a)^s \qquad (y, a \ge 0).
\]
The coefficient in the lemma is \(c_{N-1}(n-e)\) for \(n \ge e\), and both bounds, at \(n-e\), are
\(n^{N-1}/((N-1)!\,a_1 \cdots a_N) + O(n^{N-2})\).
\end{proof}

\begin{proposition}
\leanmark{\leanfull}With \(\Bm_k\) as in \eqref{eq-Ba},
\[
\Cm_{k,n} = \Bm_k\, \frac{n^{mk-1}}{(mk-1)!} + O(n^{mk-2}) .
\]
\label{prop-count}
\end{proposition}

\begin{proof}
We specialise Lemma \ref{lem-stanley} to \(x_c = x\) for the cells of the first row and \(x_c = 1\)
for the others. Setting variables equal to \(1\) is not in general allowed in a power series, so
we argue with the bijection \(\pi \mapsto (T(\pi),d)\) of that proof: \(u_t(T)\) contains
\(a_t(T)\) cells of the first row, so the first row of \(\pi\) has sum \(\sum_t d_t\, a_t(T)\).
Since \(a_t(T) \ge 1\) for all \(t\), because the entry \(1\) of \(T\) is in the first row, each
\(d_t\) is at most that sum, and the geometric series can be summed in the one variable \(x\).
With \(e(T) = \sum_{t \in D(T)} a_t(T)\), this gives
\begin{equation}
\label{eq-gf}
\leaneq{\leanfull}
\sum_{n \ge 0} \tCm_{k,n}\, x^n = \sum_{T \in \SYT(k^m)} \frac{x^{e(T)}}{\prod_{t=1}^{mk}\bigl(1-x^{a_t(T)}\bigr)} .
\end{equation}
Since \(a_1(T) = 1\), Lemma \ref{lem-coeff} applies to each term with \(N = mk\). Adding the terms,
\[
\tCm_{k,n} = \Bm_k\, \frac{n^{mk-1}}{(mk-1)!} + O(n^{mk-2}).
\]
In the same way
\(\tCm_{k-1,n} = O(n^{m(k-1)-1}) = O(n^{mk-2})\) if \(k \ge 2\), and \(\tCm_{0,n} = 0\) for
\(n \ge 1\). Now use \eqref{eq-diff}.
\end{proof}

\begin{proof}[Proof of the first part of Theorem \ref{thm-main}]
By conjugation,
\(\sum_n c_{k,n}x^n = x^k/\prod_{i=1}^{k}(1-x^i)\), so Lemma \ref{lem-coeff} gives
\(c_{k,n} = n^{k-1}/((k-1)!\,k!) + O(n^{k-2})\). Divide Proposition \ref{prop-count} by this; both
leading constants are positive.
\end{proof}

\begin{proof}[Proof of \eqref{eq-Bb}]
Let \(T \in \SYT(k^m)\). Let \(\tau_1 < \dots < \tau_k\) be the entries of its
first row and \(e_1 < \dots < e_{(m-1)k}\) the other entries. Let \(T'\) be rows \(2,\dots,m\) of
\(T\) with \(e_i\) replaced by \(i\), a tableau of shape \((k^{m-1})\), and let
\(s_i = \#\{j : \tau_j < e_i\}\), the number of first-row cells filled before \(e_i\) is entered.
The sequence \(s\) is weakly increasing with values at most \(k\). The cell \((2,j)\) lies below
\((1,j)\), so if \(e_i\) is its entry then \(s_i \ge j\). Given \(i\), let \(j = a_i(T') \ge 1\); the
cell \((2,j)\) of \(T\) holds some \(e_{i'}\) with \(i' \le i\), so \(s_i \ge s_{i'} \ge j\). Hence
\(s \in S(T')\).

Conversely, let \(T'\) and \(s \in S(T')\) be given. The relative order of the entries of \(T\) is
forced: the \(j\)-th cell of row \(1\) precedes the \(i\)-th of the other cells if and only if
\(j \le s_i\). This is a linear order because \(s\) is weakly increasing, and it is a tableau: the
only condition not inherited from \(T'\) is that \((1,j)\) precedes \((2,j)\), and \((2,j)\) is the
\(i\)-th other cell for the least \(i\) with \(a_i(T') = j\), where \(s_i \ge j\). So
\(T \mapsto (T',s)\) is a bijection from \(\SYT(k^m)\) to the pairs with \(s \in S(T')\).

When \(\tau_j\) is entered \(a_t = j\), and when \(e_i\) is entered \(a_t = s_i\). Therefore
\(\prod_t 1/a_t(T) = (1/k!)\prod_i 1/s_i\).
\end{proof}

\begin{corollary}
\leanmark{\leanfull}For \(k \ge 1\), \(A_{k,n} = G_k\, n^k\,(1+O(1/n))\) with \(G_k\) as in \eqref{eq-ar}.
\label{cor-ar}
\end{corollary}

\begin{remark}
\leanmark{\leannone}After Lemma \ref{lem-stanley}, the only property of the rectangle that is used is \(a_1(T)=1\): the
poset has a greatest element under the order of the values, and it carries the variable \(x\).
This is what makes \(x=1\) the only pole of full order. In \autocite{azam_generating_2023} the pole
order is obtained by induction through the recursion (their Proposition 17), and the
identification of the constant with the sum over \(I_k\) takes a further argument (their Propositions 18 and 21, with Lemmas 19 and 20); here the sum is read off the tableaux. The same proof gives the theorem, with the same constant \(\Gm_k\), for the average over the
partitions of \(n\) with at most \(k\) parts: \(\tCm_{k,n}\) and the number of such partitions
have the same leading terms as \(\Cm_{k,n}\) and \(c_{k,n}\). For \(m=2\) the leading term is stated in
\autocite[Remark 6]{azam_generating_2023}.
\label{rem-method}
\end{remark}
\section{Two columns}
\label{sec-two}
\begin{proof}[Proof of Proposition \ref{prop-two}]
Read \(T \in \SYT(2^m)\) as a word in \(U\) and \(D\): the
\(t\)-th letter is \(U\) if \(t\) is in the first column and \(D\) if it is in the second. This is
the usual bijection onto the Dyck paths with \(2m\) steps: the column condition says that every
prefix has at least as many \(U\) as \(D\). A \emph{return} of a Dyck path is a step \(D\) that ends on
the axis, the last step included. If the path begins with exactly \(j\) steps
\(U\), then the entry of the cell \((1,2)\) is \(j+1\), so \(a_t = 1\) for \(t \le j\) and
\(a_t = 2\) for the other \(2m-j\) values of \(t\). Hence
\[
\Bm_2 = \frac{1}{4^m}\sum_{j \ge 1} 2^j N(j),
\]
where \(N(j)\) is the number of Dyck paths with \(2m\) steps that begin with \(U^jD\). Deleting
this prefix, \(N(j)\) is the number of paths with \(2m-j-1\) steps from height \(j-1\) to height
\(0\) that stay weakly above \(0\). A Dyck path with exactly \(j\) returns to the axis is
\(UP_1D\,UP_2D \cdots UP_jD\) with Dyck paths \(P_i\), and \(P_1 D P_2 D \cdots D P_j\), started at
height \(j-1\), is a path of the kind just described; every such path arises once, by cutting it
at its first visits to the heights \(j-2,\dots,0\). So \(N(j)\) is the number of Dyck paths with
\(j\) returns, and \(\sum_j 2^j N(j)\) counts the Dyck paths in which each of the arches \(UP_iD\)
is either kept or reflected in the axis. These are all the paths with \(m\) steps \(U\) and \(m\)
steps \(D\), cut at their returns to the axis, and there are \(\binom{2m}{m}\) of them.
\end{proof}

The formal proof of Proposition \ref{prop-two} takes another route, which avoids the count of
paths by their returns.

\begin{proof}[Second proof of Proposition \ref{prop-two} (the formalised route)]
Use the words of the first proof, and give a letter the weight \(1/2\) if it is the first \(D\)
or comes after it, and the weight \(1\) otherwise; the weight of \(T\) is the product of the
weights of its letters. Call a word a \emph{ballot word} if each of its prefixes has at least as many
\(U\) as \(D\), and let \(g(n,q)\) be the total weight of the ballot words of length \(n\) with
\(q\) letters \(D\), so that \(\Bm_2 = g(2m,m)\). Removing the last letter gives, for
\(2q \le n+1\),
\[
g(n+1,0) = g(n,0), \qquad g(n+1,q) = \tfrac12\, g(n,q) + \tfrac12\, g(n,q-1) \quad (q \ge 1),
\]
while \(g(n+1,q) = 0\) for \(2q > n+1\). By induction on \(n\) and Pascal's rule,
\[
g(n,q) = \frac{1}{2^n} \sum_{i=q}^{n-q} \binom{n}{i},
\]
and for \(n=2m\), \(q=m\) the sum is the single term \(\binom{2m}{m}\).
\end{proof}
\section{An example: triples, \(k=2\)}
\label{sec-example}
\leanmark{\leannative}Let \(m=3\) and \(k=2\), so that \(C^{(3)}_{2,n}\) counts the triples
\(\nu \subseteq \mu \subseteq \lambda\) with \(\lambda \vdash n\) of length \(2\), and
\(c_{2,n} = \lfloor n/2 \rfloor\). For \(n=4\) there are \(25\) triples with \(\lambda=(3,1)\) and
\(20\) with \(\lambda=(2,2)\), so \(C^{(3)}_{2,4} = 45\) and \(A^{(3)}_{2,4} = 45/2\).

\leanmark{\leannative}The five tableaux of shape \((2^3)\), written row by row, are listed in Table \ref{tab-ex} with the
sequence \(a_1 \dots a_6\), the set \(D(T)\), the weight \(\prod_t 1/a_t\) and the pair
\((T',s)\) of the proof of \eqref{eq-Bb}.

\begin{table}[htbp]
\caption{\label{tab-ex}The tableaux of shape \((2^3)\).}
\centering
\begin{tabular}{cccccc}
\toprule
\(T\) & \(a_1 \dots a_6\) & \(D(T)\) & weight & \(T'\) & \(s\)\\
\midrule
\(14\,\vert\,25\,\vert\,36\) & \(111222\) & \(\{3\}\) & \(1/8\) & \(13\,\vert\,24\) & \(1122\)\\
\(13\,\vert\,25\,\vert\,46\) & \(112222\) & \(\{2,4\}\) & \(1/16\) & \(13\,\vert\,24\) & \(1222\)\\
\(12\,\vert\,35\,\vert\,46\) & \(122222\) & \(\{4\}\) & \(1/32\) & \(13\,\vert\,24\) & \(2222\)\\
\(13\,\vert\,24\,\vert\,56\) & \(112222\) & \(\{2\}\) & \(1/16\) & \(12\,\vert\,34\) & \(1222\)\\
\(12\,\vert\,34\,\vert\,56\) & \(122222\) & \(\emptyset\) & \(1/32\) & \(12\,\vert\,34\) & \(2222\)\\
\bottomrule
\end{tabular}
\end{table}

\leanmark{\leanfullnative}By \eqref{eq-Ba}, \(B^{(3)}_2 = 1/8 + 1/16 + 1/32 + 1/16 + 1/32 = 5/16\), in agreement with
Proposition \ref{prop-two}. By \eqref{eq-Bb}: the tableau \(T' = 13\,\vert\,24\) has \(a_1 \dots a_4 = 1122\) and
\(S(T') = \{1122, 1222, 2222\}\), contributing \(1/4+1/8+1/16 = 7/16\); the tableau
\(T' = 12\,\vert\,34\) has \(a_1 \dots a_4 = 1222\) and \(S(T') = \{1222, 2222\}\), contributing \(3/16\); and
\(\tfrac{1}{2!}(7/16+3/16) = 5/16\). Hence \(G^{(3)}_2 = \tfrac{2!\,1!}{5!}\cdot\tfrac{5}{16} = 1/192\).

\leanmark{\leannative}Here \eqref{eq-gf} and \eqref{eq-diff} give
\begin{align*}
\sum_{n} \widetilde{C}^{(3)}_{2,n} x^n &= \frac{1+2x+4x^2+2x^3+x^4}{(1-x)^6(1+x)^5},\\
\sum_{n} C^{(3)}_{2,n} x^n &= \frac{x^2\,(6+8x+x^2-6x^3-2x^4+2x^5+x^6)}{(1-x)^6(1+x)^5},
\end{align*}
with a pole of order \(6 = mk\) at \(x=1\) and of order \(5\) at \(x=-1\). The convergence is slow:
the ratio \(A^{(3)}_{2,n}/(n^4/192)\) is \(4.262\) at \(n=10\), \(1.193\) at \(n=100\) and
\(1.018\) at \(n=1000\), consistent with a relative error of order \(1/n\).
\section{Questions}
\label{sec:org1dd8269}

Is there a product or determinant formula for \(\Bm_k\) when \(k \ge 3\)? The numbers
\(B^{(3)}_k\) are \(1\), \(5/16\), \(1579/34992\), \(2072569/573308928\) for \(k = 1,\dots,4\).
The term of second order in Theorem \ref{thm-main} is not small for moderate \(n\). It could be
extracted in the same way, from the next coefficient of the expansions of \eqref{eq-gf} and of
\(\sum_n c_{k,n}x^n\) at \(x=1\). For \(k \ge 3\) no other pole matters, since the values \(1\),
\(2\) and \(3\) all occur among the \(a_t(T)\); for \(k=2\) both series also have a pole at
\(x=-1\) of order one less than at \(x=1\).
\section*{Acknowledgements}
\label{sec:orgbe36a79}
Claude (Anthropic) assisted substantially with this note; what that amounted to is set out in the
disclosure section below. The computations were done in SageMath \autocite{sagemath}.
\section*{Disclosure of AI assistance}
\label{sec:org016a3ae}
This manuscript was produced with substantial assistance from a large language model, Claude
(Anthropic): Claude Fable 5.2, in November 2026.

\begin{itemize}
\item The author posed the problem: to extend the results of \autocite{snellman_multi_2025} from pairs
to longer multichains, and then to find the analogue of Theorem 3 of
\autocite{azam_generating_2023}.
\item The model found the statement of Theorem \ref{thm-main}, with both expressions for \(\Bm_k\),
and Proposition \ref{prop-two}, and the proofs given here.
\item The model wrote the SageMath code, computed Table \ref{tab-G} and the example of Section
\ref{sec-example}, and wrote the Lean formalisation; the definitions and the statements were
fixed, and tested on small cases, before any proof was attempted.
\item The text was drafted by the model, with the author directing the exposition.
\end{itemize}

The constants in Table \ref{tab-G} were computed from \eqref{eq-Ba} and from \eqref{eq-Bb} by
separate programs, the entry \(m=4\), \(k=4\) from \eqref{eq-Ba} only, and compared with the
leading coefficient of the rational generating function at \(x=1\) for small \(k\) and \(m\). A
second run of the model, asked to referee the draft, recomputed the table and the example with
code of its own. The column \(m=2\) was compared with Table 1 of \autocite{azam_generating_2023}.
The author has verified the results to the best of his ability, made the final edits, and
assumes full responsibility for any errors, mistakes or gaps that remain.
\section*{Code and data availability}
\label{sec:org1e53a9e}
The Lean formalisation, the SageMath programs with their output, and a list of the cited
statements quoted from their sources are in the repository
\url{https://github.com/JaTaSn/young-lattice-multichains}. The directory \texttt{lean/} pins the versions of
Lean and Mathlib; its \texttt{README.md} pairs each formal statement with the result it formalises, and
\texttt{Challenge.lean} states the main results in a form that can be read without the proofs.
\section*{Funding}
\label{sec:org1a00bdf}
This research received no external funding. It was carried out as part of the author's regular
duties at Linköping University, within the fraction of that employment allocated to research.
\section*{Conflicts of interest}
\label{sec:orgfcf76f5}
The author declares no financial or non-financial conflicts of interest. The nature and extent of
large-language-model assistance in producing this manuscript is disclosed in full in the section
above.

\printbibliography
\end{document}